\documentclass[10pt]{article}   %% Standard LaTeX.
\usepackage[letterpaper, margin=1 in]{geometry}
\usepackage{amsmath}
\usepackage{mathtools}
\usepackage{amssymb}
\usepackage{bm}      
\usepackage{amsthm}
\usepackage{enumerate}
\usepackage{tikz}
\usepackage{tikz-cd}
\usepackage{graphicx}        %% Enable for eps figures
\usepackage{verbatim}
\usepackage{comment}
\usepackage{quiver}
\usepackage{setspace}
\usepackage[all,cmtip]{xy}
\usepackage{comment}
\usepackage[backend=biber, natbib=true, style=numeric]{biblatex}

\theoremstyle{plain}
\newtheorem{thm}{Theorem}[section]
\newtheorem{lm}[thm]{Lemma}

\theoremstyle{definition}

\newtheorem{define}[thm]{Definition}

\newcommand{\R}{\mathbb{R}}
\newcommand{\Q}{\mathbb{Q}}
\newcommand{\Z}{\mathbb{Z}}

\newcommand{\B}{\mathcal{B}}
\newcommand{\M}{\mathcal{M}}

\DeclareMathOperator{\Ext}{Ext}
\DeclareMathOperator{\Mono}{Mono}
\DeclareMathOperator{\planes}{planes}

\newcommand{\id}{\mathrm{id}}

\newcommand{\inv}{^{-1}}

\newcommand{\Hom}{\operatorname{Hom}}

\title{Immersibility of manifolds is decidable in odd codimension}

\author{Helen Epelbaum}

\begin{document}

\emergencystretch 3em

\maketitle

\begin{abstract}
  Given a smooth map $f:M\rightarrow N$ of closed oriented smooth manifolds, is there an immersion homotopic to $f$? We provide an algorithm that decides this when the codimension of the  manifolds is odd.
\end{abstract}

\section{Introduction}

Given a pair of smooth manifolds, when can we immerse one in the other? A lot is known in the case of immersibility into $\R^{n}$. The Whitney immersion theorem tells us that any manifold of dimension $m$ can be immersed in $R^{2m-1}$. In 1985, Cohen strengthened this result proving the immersion conjecture: any manifold of dimension $m$ can be immersed in $\R^{2m-\alpha(m)}$ where $\alpha(m)$ is the number of $1$s in the binary expansion of $m$ \cite{cohen1985immersion}. This bound is tight: for any $m$ there are manifolds of dimension $m$ that cannot be immersed in dimension $2m-\alpha{m}$. In smaller codimension, we want to ask the question about whether immersibility is {\it decidable}, that is, for which pairs $(m,n)$ can a computer algorithm always determine whether a manifold of dimension $m$ immerses in $\R^{n}$. This was studied in \cite{manin2023algorithmic} which investigates both smooth and PL manifolds. For smooth manifolds, they prove among other things that immersibility of an $m$-manifold into $\R^{n}$ is decidable when $n-m$ is odd, and it is this result we will generalize here.

More generally, we might be interested in the question of when $M$ can be immersed in $N$, for arbitrary smooth manifolds $M,N$. Here we might hope to compute the set of immersions, up to regular homotopy (homotopy through immersions.) Here we will end up with difficulties arising from the difficulty of studying homotopy classes of maps between manifolds in general. To control this we will discuss a modification of the problem in which we are looking for immersions with prescribed homotopy behavior. That is, we will consider the following decision problem: given a pair of oriented closed smooth manifolds $(M,N)$ and a smooth map $f:M\rightarrow N$, is there an immersion homotopic to $f$? The main result of this paper will be a generalization of the above result for this problem. In particular we prove the following theorem:

\begin{thm}\label{mainresult}
    There is an algorithm that on input $(M, N, f)$, where $M,N$ are smooth oriented manifolds of dimension $m,n$ respectively and $f:M\rightarrow N$ is a smooth map between them, decides whether there is an immersion $g:M\rightarrow N$ such that $g\simeq f$, as long as $n-m$ is odd.
\end{thm}

This proof will have several steps. First, we use the $h$-principle of Hirsch and Smale to reduce the question to a homotopy-theoretic lifting problem.

To decide on the existence of a lift then, we will use some tools from rational homotopy theory. The idea here is that if we ignore finite homotopy group obstructions, we can put an algebraic structure on the set of possible lifts on each stage of the relevant Moore-Postnikov tower, and this allows us to construct a lift, if one exists, using obstruction theory. The algorithm here is very much like that in \cite{manin2022rational}.

In section \ref{effectiverep} we will review effective representation of smooth manifolds, so it is clear how to input such an object into an algorithm. In section \ref{hprinciple} we will show how to convert the immersion question into a homotopy lifting problem. Section \ref{ratfibhmspacetheory} will be a review of the relevant algebraic tools for the lifting algorithm, which will be presented in detail in section \ref{liftingalgsection}. We then present a full proof of theorem \ref{mainresult} in section \ref{mainresultsection}.

\section{Effective Representation}\label{effectiverep}

We will quickly summarize the needed results here, the details of which can be found in \cite{manin2023algorithmic}.

\begin{thm}\label{comphtpyresults}
    \begin{enumerate}[(a)]
        \item There is an algorithm which on input a simplicial complex $X$ can compute generators and relations for $\pi_k(X)$, as well as simplicial representatives for each generator
        \item There is an algorithm which on input a map of simply connected finite simplicial complexes $Y\rightarrow B$ computes the relative Moore-Postnikov tower to any finite stage, as well as the cohomology of each stage, and the maps of cohomology induced by each $P_n\rightarrow P_{n-1}$.
        \item There is an algorithm such that given a diagram
        \[\xymatrix{
            A \ar[r] \ar@{^(->}[d] & P_n \ar@{->>}[d] \\
            X \ar[r] \ar@{-->}[ru] & P_{n-1},
        }\] where $(X, A)$ is a finite simplicial pair and $P_n\twoheadrightarrow P_{n-1}$ is a Moore-Postnikov stage, computes the obstruction to filling in the dotted arrow (in $H^n(X, A; \pi_n(P^n))$) and if the obstruction vanishes, constructs a lifting extension. 
    \end{enumerate}
\end{thm}

We will also need to specify how we will model smooth manifolds algorithmically. In general, it is not decidable whether a given $n$-dimensional simplicial complex is homeomorphic to a smooth manifold. Again following \cite{manin2023algorithmic}, we will input a manifold as a finite simplicial complex together with a choice of polynomial map with rational coefficients to $\R^{N}$ for each top dimensional simplex, fixing some large $N$, and such that the derivatives of each map are nonsingular and agree on the boundaries of adjacent simplices. This gives a $C^{1}$ triangulation, and since the category of $C^{1}$-manifolds is equivalent to the category of $C^{\infty}$-manifolds, and $C^{1}$ immersions can be approximated by $C^{\infty}$ ones, this is sufficient. We also note that given such a collection of data, whether or not it represents a manifold is decidable.

We will need one more theorem from \cite{manin2023algorithmic} which will help convert our geometric problem into a homotopy theoretic one.
\begin{thm}\label{compgeoresults}
    Given a manifold $M$ as above, there is an algorithm that computes the classifying map of the tangent bundle. In particular, it is possible to compute a simplicial complex structure on $BSO(n)$, and a simplicial approximation of the classifying map.
\end{thm}

\section{Converting to a Lifting Problem}\label{hprinciple}

We now turn to the problem of reducing the question of immersibility to a homotopy lifting problem. The first step is to use the h-principle of Hirsch and Smale\textemdash{}the existence of an immersion homotopic to $f:M\rightarrow N$ is equivalent to the existence of a tangent bundle monomorphism $F:TM\rightarrow TN$ sitting over $f:M\rightarrow N$ \cite{gromov2013partial}.

To convert this to a lifting problem, we will construct a bundle, which we will denote by $\pi:\Mono(TM, TN)\twoheadrightarrow M\times N$, for any given manifolds $M, N$ of dimension $m, n$. The fibers of this bundle will be homotopic to the real Stiefel manifolds $V(m, n)$, the space of orthonormal $m$-frames in $\R^n$. Over each $(p, q)$ we want to think of this as the space of immersions (we use here heavily the fact that $GL(n)$ deformation retracts onto $O(n)$ to speak of an orthogonal structure and simplify some computation) of $T_p(M)$ into $T_q(N)$. To construct this space explicitly we will describe a system of transition maps on local trivializations.

Here we have to set up some notation. We will consider the tangent bundle over $M$ constructed as a collection of charts with transition maps, in particular a collection of opens $\mathcal{U}$ over $M$ with transition maps $\varphi_{ij}$ on $(U_i\cap U_j)\times V(m, n)\simeq \pi^{-1}(U_{i}\cap U_{j})$ satisfying that $\varphi_{jk}\circ\varphi_{ij}=\varphi_{ik}$ on $\pi^{-1}(U_i\cap U_j\cap U_k)$. Similarly we have the tangent bundle over $N$, with local trivializations $\mathcal{V}$ and transition maps $\psi_{ij}$. Then we construct a $V(m, n)$-bundle on $M\times N$ as follows: we construct an open cover $\mathcal{W}$ on $M\times N$ by taking the open sets $W=U\times V$ for each $U, V$ in $\mathcal{U},\mathcal{V}$ respectively. An intersection $W_i\cap W_j$ can be written as $(U_{i_1}\cap U_{j_1})\times (V_{i_2}\cap V_{j_2})$ and then the transition map can be defined on orthogonal frames $O$ by $\zeta_{i,j}(O)=\psi_{i_2 j_2}^*O\varphi_{i_1, j_1}$, which provides a linearly independent frame. A straightforward calculation shows that for any triple intersection we have $\zeta_{jk}\zeta_{ij}=\zeta_{ik}$. It is here that we use that we can, for example by singular value decomposition, retract $GL(n)$ onto $O(n)$ to consistently create orthonormal frames from linearly independent ones. Then we have constructed an open cover on $M\times N$ and described a coherent system of transition maps between the trivialization at each open, and hence we have a $V(m, n)$ bundle on $M\times N$. It remains to show that it correctly parametrizes the tangent bundle monomorphisms.

\begin{lm}\label{liftsaremonos}
    Fix a smooth map $f:M\rightarrow N$. The set of homotopy classes of orthogonal tangent bundle monomorphisms $TM\rightarrow TN$ over $f$ is in bijective correspondence with homotopy classes of lifts of the triangle:
    % https://q.uiver.app/#q=WzAsMyxbMCwxLCJNIl0sWzEsMSwiTVxcdGltZXMgTiJdLFsxLDAsIk1vbm8oVE0sIFROKSJdLFswLDEsImlkXFx0aW1lcyBmIl0sWzIsMV0sWzAsMiwiIiwyLHsic3R5bGUiOnsiYm9keSI6eyJuYW1lIjoiZGFzaGVkIn19fV1d
\[\begin{tikzcd}
	& {\Mono(TM, TN)} \\
	M & {M\times N}
	\arrow["{id\times f}", from=2-1, to=2-2]
	\arrow[from=1-2, to=2-2]
	\arrow[dashed, from=2-1, to=1-2]
\end{tikzcd}\]
\end{lm}

\begin{proof}
    Fix open covers $\mathcal{U}$ and $\mathcal{V}$ of $M$ and $N$ admitting orthonormal frames.
    
    Let $[\phi]$ be a homotopy class of tangent bundle monomorphism containing a specific monomorphism $\phi$. Then pick some $\phi\in [\phi]$ and we will construct a lift $\psi: M\rightarrow \Mono(TM, TN)$.  Then for each $u\in\mathcal{U}$ and each $v\in \mathcal{V}$ with $f(u)\cap v$ nonempty $\phi$ provides, for each $p\in u$ with $f(p)\in v$, an $m$-frame $O$ over $(p, f(p))$ viewed as a point in the local trivialization over $u\times v$, and again using a retraction to $O(n)$ we can make the frame orthonormal. If we had some $u'$ such that $p\in u'$ and $v'$ such that $f(p)\in v'$ then we would produce a frame $O'$, and these would be related precisely by $O'=\psi^* O \varphi$ where $\varphi$ is the transition map from $u$ to $u'$ for the tangent bundle on $M$ and $\psi$ is the transition map from $v$ to $v'$ for the tangent bundle on $N$. This is exactly the transition map from $u\times v$ to $u'\times v'$ for $\Mono(TM, TN)$ on $M\times N$, so this produces a lift. If we had picked a different $\phi'\in[\phi]$ they would be related by a homotopy, and again reducing to coordinates over $\mathcal{U}\times \mathcal{V}$, we could produce a homotopy of lifts.

    Conversely, given a homotopy class $[g]$ of lift of the above diagram, we pick some lift $g$ in the class and produce a tangent bundle monomorphism. Again, for each $p$ in $M$, each $u\in U$ containing $p$ and each $v\in \mathcal{V}$ containing $f(p)$ the lift $g$ provides an $m$-frame over $(p, f(p))$ which we can view as a map from $T_pM$ to $T_{f(p)}N$ in the coordinates of $u$ and $v$. Again that this is a coherent choice of frame across all choices of $u,v$ follows from the fact that the transition maps agree. Finally, a different $g' \in [g]$ would be related by a homotopy to $g$, and this can again be written out in coordinates.
\end{proof}

We want to take this a step further however, as we would like to produce a lifting problem where the fibration is uniform over all $N$, once dimensions are fixed. To do this, we construct a bundle, in the same way as $\Mono(TM, TN)$ over $BSO(m)\times BSO(n)$, which we will denote $p_{m,n}:\Mono(m-\planes, n-\planes)\twoheadrightarrow BSO(m)\times BSO(n)$. In particular, we can again start with the universal bundles $ESO(m)\twoheadrightarrow BSO(m)$ and $ESO(n)\twoheadrightarrow BSO(n)$, written out in coordinates over some systems of local trivializations. We prove the following theorem.

\begin{lm}\label{bundlepullback}
    There is a bundle $p_{n.m}:\Mono(m-\planes, n-\planes)\twoheadrightarrow M\times N$ such that for all smooth manifolds $M, N$ of dimension $m, n$ respectively, there is a commutative diagram:
    % https://q.uiver.app/#q=WzAsNCxbMCwwLCJNb25vKFRNLCBUTikiXSxbMSwwLCJNb25vKG0tcGxhbmVzLCBuLXBsYW5lcykiXSxbMCwxLCJNXFx0aW1lcyBOIl0sWzEsMSwiQlNPKG0pXFx0aW1lcyBCU08obikiXSxbMCwxXSxbMCwyLCIiLDIseyJzdHlsZSI6eyJoZWFkIjp7Im5hbWUiOiJlcGkifX19XSxbMSwzLCJwX3ttLG59IiwwLHsic3R5bGUiOnsiaGVhZCI6eyJuYW1lIjoiZXBpIn19fV0sWzIsMywiXFxrYXBwYV9tXFx0aW1lc1xca2FwcGFfbiIsMl1d
\[\begin{tikzcd}
	{Mono(TM, TN)} & {\Mono(m-\planes, n-\planes)} \\
	{M\times N} & {BSO(m)\times BSO(n)}
	\arrow[from=1-1, to=1-2]
	\arrow[two heads, from=1-1, to=2-1]
	\arrow["{p_{m,n}}", two heads, from=1-2, to=2-2]
	\arrow["{\kappa_m\times\kappa_n}"', from=2-1, to=2-2]
\end{tikzcd}\]
where $\kappa_m, \kappa_n$ are the classifying maps for the tangent bundles of $M, N$ respectively. Furthermore, the above diagram is a pullback square.
\end{lm}

\begin{proof}
    Again we pick a covering collection of opens $\mathcal{U}$ on $BSO(m)$ and $\mathcal{V}$ on $BSO(n)$ admitting local trivializations, and consider the open cover $\mathcal{W}$ of $BSO(m)\times BSO(n)$ by taking sets of the form $W=U\times V$. Then over each of these sets we have a trivial $V(m, n)$ bundle and we construct transition maps from $U\times V$ to $U'\times V'$ as before: if $\varphi$ is the transition map from $U$ to $U'$ and $\psi$ is the transition map from $V$ to $V'$ then over each point in $(U\times V)\cap (U'\times V')=(U\times U')\cap (V\times V')$ we have the transition map $\zeta(O)=\psi^* O\varphi$.

    The map $\Mono(TM, TN)\rightarrow \Mono(m-\planes, n-\planes)$ works as follows: we can pull back the collection $\mathcal{U}$ on $BSO(m)$ across $\kappa_m$ to form a collection of opens on $M$, and the local trivializations of $ESO(m)\twoheadrightarrow BSO(m)$ pull back to trivializations of the tangent bundle of $M$. Similarly, the collection of opens $\mathcal{V}$ pulls back across $\kappa_n$ along with the trivializations of $ESO(n)\twoheadrightarrow BSO(n)$. Then for any point $(p, q)\in M\times N$ we can pick some $\kappa_m\inv(u)\times \kappa_n\inv(v)$ containing $(p, q)$ and we can use this choice of trivialization to map the fiber over $(p, q)$ to the fiber over $\kappa_m(p)\times \kappa_n(q)$. That this describes a map coherently over all choices of trivialization of $BSO(m)$ and $BSO(n)$ follows from the fact that the tangent bundles of $M$ and $N$ are pullbacks of $ESO(m)\twoheadrightarrow BSO(m)$ and $ESO(n)\twoheadrightarrow BSO(n)$ respectively.
    
    Finally, this is a pullback square because firstly it would factor through the pullback, and the map to the pullback would be a map of fiber bundles with the same base, and hence a projection of fibers, and since the fiber is the same it is hence a homeomorphism. 
\end{proof}

We can put this together with lemma \ref{liftsaremonos} to get the following theorem:

\begin{thm}\label{universalliftsaremonos}
    Given $M, N, f$ where $M$ is a smooth $m$-dimensional manifold, $N$ is a smooth $n$-dimensional manifold and $f$ is a smooth map between them, with $p_{m,n}:\Mono(m-\planes, n-\planes)\twoheadrightarrow BSO(m)\times BSO(n)$ as in lemma \ref{bundlepullback}, then there is a bijective correspondence between homotopy classes of lifts of the diagram:
    % https://q.uiver.app/#q=WzAsMyxbMCwxLCJNIl0sWzEsMSwiQlNPKG0pXFx0aW1lcyBCU08obikiXSxbMSwwLCJNb25vKG0tcGxhbmVzLCBuLXBsYW5lcykiXSxbMCwxLCIoXFxrYXBwYV9tXFx0aW1lc1xca2FwcGFfbilcXGNpcmMoaWRcXHRpbWVzIGYpIiwyXSxbMiwxLCJwX3ttLG59IiwwLHsic3R5bGUiOnsiaGVhZCI6eyJuYW1lIjoiZXBpIn19fV0sWzAsMiwiIiwyLHsic3R5bGUiOnsiYm9keSI6eyJuYW1lIjoiZGFzaGVkIn19fV1d
\[\begin{tikzcd}
	& {Mono(m-planes, n-planes)} \\
	M & {BSO(m)\times BSO(n)}
	\arrow["{(\kappa_m\times\kappa_n)\circ(id\times f)}"', from=2-1, to=2-2]
	\arrow["{p_{m,n}}", two heads, from=1-2, to=2-2]
	\arrow[dashed, from=2-1, to=1-2]
\end{tikzcd}\]
and homotopy classes of tangent bundle monomorphisms over $f$.
\end{thm}

Then we have successfully converted the question of whether there is an immersion homotopic to a given map to the question of whether there is a lift of a diagram. To find such a lift, we will need some tools from rational homotopy theory. We turn to that in the next section.

\section{Rational Fibrewise HM-Spaces}\label{ratfibhmspacetheory}

The main tools we will need from rational homotopy theory are the notion of a minimal model, and its relative version,  developed by Sullivan in \cite{sullivan1977infinitesimal}. We will give a brief summary here, for a more thorough introduction, see \cite{felix2012rational}. The key idea is that given a simply connected CW complex $X$, we can construct a rational graded commutative differential algebra in a way that institutes a duality between the category of simply connected CW complexes (up to rational equivalence) and the category of $1$-connected rational cdgas (up to quasi-isomorphism.) In particular, we construct the {\it minimal model} $\mathcal{M}_X$ as follows. We begin with a graded vector space $W$ whose $i^{(th)}$ graded piece $W^{(i)}$ is given by $(\pi_i(X)\otimes \Q)^*$. Note that since we started with a simply connected space, the resulting vector space is trivial below dimension $2$, and the resulting dga is said to be {\it $1$-connected} or {\it simply connected.} Then as an algebra, the minimal model of $X$ is $\wedge W$, the free graded commutative algebra generated by $W$. Then a differential on $\mathcal{M}_X$ is determined by maps $d_i:W^{(i)}\rightarrow \mathcal{M}_X^{(i+1)}$, which can be constructed by dualizing the $k$-invariants of the Postnikov tower, after tensoring with $\Q$. Crucially this will result in maps which land in the subalgebra generated by elements of $W^{(k)}$ for $k<i$, and so for each $w\in W^{(i)}$ we have $d(w)\in \wedge^{\geq 2} \oplus_{j=2}^{i-1}W^{(j)}$. Such a dga is said to be {\it minimal}.

A space is modelled by a minimal model, and a fibration is modelled by a {\it relative minimal model.} Suppose we have a fibration of simply connected CW-complexes $p:Y\twoheadrightarrow B$, with fiber $F$ which is also simply connected. Then we can construct a {\it relative minimal model.} This consists of the following:
\begin{itemize}
    \item $\mathcal{M}_B=(\wedge W_B, d_B)$ a minimal model for the base space $B$
    \item A graded vector space $W_F$ where each $i^{th}$ graded piece is given by $(\pi_i(F)\otimes \Q)^*$
    \item A differential $d$ on $\mathcal{M}_X\otimes_\Q \wedge W_F$ which restricts to $d_B$ under the standard inclusion $\mathcal{M}_X\rightarrow \mathcal{M}_X\otimes_\Q \wedge W_F$ given by $x\mapsto x\otimes 1$.
\end{itemize}
Again this should be {\it minimal} which amounts to that $d(W^{(i)})\subset \mathcal{M}_X\otimes_\Q \wedge \oplus_{j=2}^{i-1}W_F^{(j)}$, and for this to be a model of the fibration the cohomology of the dga $(\mathcal{M}_X\otimes \wedge W_F, d)$ should be the cohomology of $Y$. Similar to the minimal model, the differential can be constructed by dualizing the Moore-Postnikov tower.

We prove a few lemmas.

\begin{lm}\label{lineq}
    Let $(A, d)$ be a minimal rational finitely generated simply connected dga. Suppose we are given a set of $k$ equations
    \begin{align*}
        dx_i=a_i+\sum_{j=1}^k a^i_j x_j
    \end{align*}
    with each $a_i$ and $a^i_j\in A$
    and a prescribed sequence $2\leq n_1\leq ... \leq n_k$. Then there is an algorithm to determine the (possibly empty) affine space of solutions with the condition $|x_i|=n_i$, presented as a solution $s=(s_1,...,s_k)$ and a basis of a space $W$ such that the affine space $s+W$ is the set of solutions in the $\Q$-vector space $V=\oplus_{i=1}^k A^{(n_i)}$.
\end{lm}

\begin{proof}
    This is simply a matter of solving a matrix equation. Indeed, with $V$ as given in the statement of the lemma, we construct a pair of maps $D,T:V\rightarrow \oplus_{i=1}^k A$ as follows: on $v=(v_1,\ldots,v_k)$ define
    \begin{align*}
        Dv&=(dv_1,...,dv_k) \\
        Tv&=(\sum_{i=1}^k a^1_i v_i,...,\sum_{i=1}^k a^k_iv_i)
    \end{align*}
    and set
    $$C=(a_1,...,a_k)$$
    so that the system of equations is equivalent to
    $$Dv=Tv+C$$
    which we can now turn into a matrix equation. In particular, as each $A^{(n_i)}$ can be represented as a rational vector space, as can $A$ itself, and since the differential is $\Q$-linear, so is $D$, as is $T$ by construction. Then picking a basis for $V$ and $A$ as $\Q$-vector spaces, we have the above equation as a matrix equation, and the space of solutions can be constructed via row reduction.
\end{proof}

Given a square 
\begin{center}
    % https://tikzcd.yichuanshen.de/#N4Igdg9gJgpgziAXAbVABwnAlgFyxMJZABgBpiBdUkANwEMAbAVxiRAB12BbOnACwDGjYAFkAvgH0AgiDGl0mXPkIoyARiq1GLNpx78hDUZIAas+SAzY8BImvKb6zVog7deg4eIkBNcwutlO1INaicdVz0PQ2MJACFZTRgoAHN4IlAAMwAnCC4kAGZqHAgkNTks3PzEexASpDIQPhg6KDYcAHcIZtaECpAcvMLi0sQAJn7B6rGRhrEKMSA
    \begin{tikzcd}
    \mathcal{M}_A                      & \mathcal{M}_Y \arrow[l]           \\
    \mathcal{M}_X \arrow[u] & \mathcal{M}_B \arrow[l] \arrow[u]
    \end{tikzcd}
    \end{center}
and a relative minimal model $(\M_B \otimes \wedge W, d)$ with a differential d which is linear through dimension $n$, where $n$ is the cohomological dimension of a CW pair $(X, A)$ for which $\M_X$ and $\M_A$ are minimal models, and $\M_X\rightarrow \M_A$ models the inclusion $A\hookrightarrow X$. 

Taking the pushout of bottom right triangle, and replacing $\M_Y$ with the given relative model we arrive at the following triangle, for which we would like to construct the dashed line:
\begin{center}
 % https://tikzcd.yichuanshen.de/#N4Igdg9gJgpgziAXAbVABwnAlgFyxMJZARgBoAGAXVJADcBDAGwFcYkQAdDgWQH0ANLhDwBbeAAIuAdxhQA5jHEB1EAF9S6TLnyEU5UsWp0mrdlz4BBNRpAZseAkTKGaDFm0SceAtUdkKEFFAAMwAnCBEkfRAcCCQyY3d2YIA9ACoQGkZ6ACMYRgAFLQddEFCsOQALHGsQ8MjEaNikACZXEw8wZkZGLNz8ovsddkYYYJqsrDAPECh6OErZTJA8sCgkAFoAZnJ1OojWmmbG9qTPLjRKrHTl7LzC4uHPcqqavZAwg8Q2mLjEBLcpk81wyfXug20jk8o3GvlUQA
\begin{tikzcd}
    & \M_X\otimes \wedge W \arrow[ld, "f^*"'] \arrow[d, dashed, bend right] \\
\M_A & \M_X \arrow[u, "\phi^*"'] \arrow[l, "i^*"]                           
\end{tikzcd}
\end{center}

\begin{lm}\label{dgalift}
    Given a diagram as above, there is an algorithm to determine whether such a dashed line exists.
\end{lm}

\begin{proof}
Fix a basis for $W$, written as a choice of elements $w^{(i)}_j$ where $i$ ranges over the positive dimensional degrees, and $j$ ranges from $1$ to the dimension of $W^{(i)}$. Constructing a map for the dashed line as in the above diagram is simply a matter of picking a target for each element $w_j^{(i)}$, as $\M_X \otimes \wedge W$ is the free commutative graded $\M_X$-algebra. The only conditions then to check on such a map $\psi$ are $d_{\M_X}\circ\psi=\psi\circ d$ and $f^*=i^*\circ \psi$. Then for each $w^{(i)}_j$ we have the equations $d_{\M_X}(\psi(w_j^{(i)}))=\psi(d(w^{(i)}_j))$ and $f^*(w^{(i)}_j)=i^*(\psi(w^{(i)}_j))$.

By assumption, $d$ is linear through the cohomological dimension of $(X, A)$, and in higher dimensions there is a unique lift of any element in $\M_A$. Then we have a fixed $\psi(w_j^{(i)})$ for $i>n$. 

Otherwise we know that $dw^{(i)}_j$ takes the form $$m^{i}_j+\sum_{k<i, 1\leq k \leq dim(W^{(i)})} m^{i,k,l}_j w^{(k)}_l$$ where each $m^{i}_j$ and $m^{i,k,l}_j\in\M_X$. Putting all of this together then, we obtain the set of equations:
\begin{align*}
    d_{\M_X}(\psi(w^{(i)}_j))=m^i_j+\sum_{k<i,1\leq k\leq dim(W^{(i)})} m_j^{i, k, l}\psi(w^{(k)}_l)
\end{align*}
and by lemma \ref{lineq} we can construct the affine space of solutions to these equations as a subspace of $V=\oplus_{i,j} A^{(i)}$ (note that we have a copy of $A^{(i)}$ for each generator in $W^{(i)}$.) We can call this space $S$.

We consider the affine space $\tilde{V}\subset V$ given by:
$$\oplus_{i,j} i^{*-1}(f^*(w^{(i)}_j))$$
and we simply have to compute $S\cap\tilde{V}$, which is the intersection of two affine subspaces of a $\Q$-vector space, which can be computed.

\end{proof}

Before we are able to prove the lifting result, we will need two different generalizations of the notion of an $H$-space.

\begin{define}
    A fibrewise $H$-space is a fibration $p:Y\twoheadrightarrow B$ with a section $e:B\rightarrow Y$ and a multiplication map $m:Y\times_B Y\rightarrow B$ which is associative up to fibrewise homotopy, and for which the section acts as an identity, that is the maps $m\circ (\id \times (e\circ p))\circ \Delta$ and $m\circ ((e\circ p)\times \id)\circ \Delta$ are fibrewise homotopic to the identity on $Y$, where $\Delta:Y\rightarrow Y\times_B Y$ is the diagonal map.
\end{define}

We note that this is stronger than simply requiring an $H$-space structure on each fiber, the fibers need to have an $H$-space in some strong uniform sense. We want to weaken this definition slightly however, to include spaces which may not have a section. A motivating example here is the Hopf fibrations. The lifting algorithm we construct will allow us to compute the existence of lifts across the fibration $p:S^7\twoheadrightarrow S^4$ with fiber $S^3$, despite the fact that this doesn't admit a fibrewise $H$-space structure.
One way of looking at a `group without identity,' is to look at the notion of a heap.

\begin{define}[Heap]
  A {\it Mal'cev operation} on a set $H$ is a ternary operation $\tau$ which satisfies
  \begin{itemize}
    \item $\forall a,b,c,d,e\in H, \tau(\tau(a,b,c),d,e)=\tau(a,\tau(b,c,d),e)=\tau(a,b,\tau(c,d,e))$
    \item $\forall a,b\in H, \tau(a,a,b)=b=\tau(b,a,a)$
  \end{itemize}
The first of these is a kind of associativity, and the second is often referred to as the Mal'cev condition. A {\it Heap} is a set with a Mal'cev operation.
\end{define}

If $(H,\tau)$ is a heap, it is a straightforward exercise to check that for any $e\in H$, the operation $a*b=\tau(a,e,b)$ turns $H$ into a group. On the other hand if we have a group $(G,*)$ the operation $\tau(a,b,c)=a*b\inv *c$ is a Mal'cev operation which turns the group into a heap. It is natural then to consider a heap to be a group with a forgotten identity element. The choice of any distinguished element recovers the structure of a group.

If an $H$-space is a `grouplike' space, we want to look at something `heaplike.' This motivates the next definition.

\begin{define}
    A fibrewise $HM$-space (for Hopf-Mal'cev) is a fibration $p:Y\twoheadrightarrow B$ together with a fibrewise homotopy Mal'cev operation, i.e. a map $\tau:Y\times_B Y \times_B Y \rightarrow Y$ for which the following diagrams are commutative, up to fibrewise homotopy:

    % https://q.uiver.app/#q=WzAsNCxbMCwwLCJZXFx0aW1lc19CWVxcdGltZXNfQiBZXFx0aW1lc19CIFlcXHRpbWVzX0JZIl0sWzEsMCwiWVxcdGltZXNfQiBZXFx0aW1lc19CIFkiXSxbMCwxLCJZXFx0aW1lc19CIFlcXHRpbWVzX0IgWSJdLFsxLDEsIlkiXSxbMCwxLCJcXHRhdVxcdGltZXMgaWRcXHRpbWVzIGlkIl0sWzAsMiwiaWRcXHRpbWVzIGlkXFx0aW1lcyBcXHRhdSIsMl0sWzIsMywiXFx0YXUiLDJdLFsxLDMsIlxcdGF1Il1d
\[\begin{tikzcd}
	{Y\times_BY\times_B Y\times_B Y\times_BY} & {Y\times_B Y\times_B Y} \\
	{Y\times_B Y\times_B Y} & Y
	\arrow["{\tau\times id\times id}", from=1-1, to=1-2]
	\arrow["{id\times id\times \tau}"', from=1-1, to=2-1]
	\arrow["\tau"', from=2-1, to=2-2]
	\arrow["\tau", from=1-2, to=2-2]
\end{tikzcd}\]

% https://q.uiver.app/#q=WzAsNCxbMCwwLCJZXFx0aW1lc19CWSJdLFsyLDAsIllcXHRpbWVzX0JZXFx0aW1lc19CWSJdLFsyLDEsIlkiXSxbNCwwLCJZXFx0aW1lc19CIFkiXSxbMCwxLCJpZFxcdGltZXMgXFxEZWx0YSJdLFsxLDIsIlxcdGF1Il0sWzAsMiwiXFxwaV8xIiwyXSxbMywxLCJcXERlbHRhXFx0aW1lcyBpZCIsMl0sWzMsMiwiXFxwaV8yIl1d
\[\begin{tikzcd}
	{Y\times_BY} && {Y\times_BY\times_BY} && {Y\times_B Y} \\
	&& Y
	\arrow["{id\times \Delta}", from=1-1, to=1-3]
	\arrow["\tau", from=1-3, to=2-3]
	\arrow["{\pi_1}"', from=1-1, to=2-3]
	\arrow["{\Delta\times id}"', from=1-5, to=1-3]
	\arrow["{\pi_2}", from=1-5, to=2-3]
\end{tikzcd}\]
where $\pi_i$ denotes projection onto the $i^{(th)}$ coordinate.
These are exactly the same axioms as for a heap of sets, we simply require only that they commute up to fibrewise homotopy rather than on the nose.

\end{define}

On sets, the choice of an identity element turns a heap into a group. An analogous result holds here.

\begin{lm}
    Let $p:Y\twoheadrightarrow B$, $\tau$ be a fibrewise HM-space, and let $e:B\rightarrow Y$ be a section of $p$. Then with multiplication map $m:Y\times_B Y\rightarrow Y$ given by $\tau \circ (\id\times (e \circ p) \times \id )\circ(id\times \Delta)$, $p$ is a fibrewise $H$-space.
\end{lm}

\begin{proof}
    We begin by observing that since $m$ is defined on elements of the fibrewise product $(a,b)\in Y\times_B Y$, we have that $p(a)=p(b)$ so that we have $(\id\times (e \circ p) \times \id )\circ(\id\times \Delta)=(\id\times (e \circ p) \times \id )\circ(\Delta \times \id)$. Then checking associativity of multiplication is a straightforward calculation. Indeed we have:
    \begin{align*}
        m\circ(m\times \id)&= \\
        \tau \circ (\id\times (e \circ p) \times \id )\circ(\id\times \Delta) \circ ((\tau \circ (\id\times (e \circ p) \times \id )\circ(\id\times \Delta))\times \id)&=\\
        \tau\circ (\tau \times \id \times \id)\circ (\id \times (e\circ p)\times \id \times (e \circ p)\times \id)\circ (\id \times \Delta \times \Delta)
    \end{align*}
    simply by substitution, and then applying the commutativity of operations on different copies of the product.
    Since the fibrewise Mal'cev operation has to also be homotopy associative, this is homotopic to
    \[
        \tau\circ (\id \times \id \times \tau)\circ (\id \times (e\circ p)\times \id \times (e \circ p)\times \id)\circ (\id \times \Delta \times \Delta)
    \]
    Applying the equality above we have that this is equal to
    \[
        \tau\circ (\id \times \id \times \tau)\circ (\id \times (e\circ p)\times \id \times (e \circ p)\times \id)\circ (\Delta \times \id \times \Delta)
    \]
    and finally this can be rewritten as
    \[
        \tau\circ (\id\times (e\circ p)\times \id)\circ (\Delta \times \id)\circ (\id \times m)    
    \]
    and again applying the above equality we have this is equal to
    \[
        \tau\circ (\id\times (e\circ p)\times \id)\circ (\id \times \Delta)\circ (\id \times m)=m\circ (\id \times m)
    \]

    To check that the section acts as an identity, we first observe that $((e\circ p) \times \id)\circ \Delta \circ (e\circ p)= \Delta \circ (e\circ p) = (\id \times (e \circ p))\circ \Delta \circ (e\circ p)$ since $(e\circ p)\circ (e\circ p)=(e\circ p)$ and for any map $f$ we have $(f\times f)\circ \Delta = \Delta \circ f$. 

    Now we compute
    \begin{align*}
        m\circ (\id \times (e\circ p))\circ \Delta&=\\
        \tau \circ (\id\times (e \circ p) \times \id )\circ(id\times \Delta) \circ  (\id \times (e\circ p))\circ \Delta &=\\
        \tau \circ (id\times \Delta) \circ (\id \times (e \circ p))\times \Delta
    \end{align*}
    and by the axioms of a fibrewise $HM$-space, $\tau\circ (\id \times \Delta)$ is fibrewise homotopic to $\pi_1$ and the above map is fibrewise homotopic to 
    \[
        \pi_1\circ (\id\times (e\circ p))\circ \Delta=\id
    \]
    An analogous argument shows that $m\circ ((e\circ p)\times \id)\circ \Delta$ is also fibrewise homotopic to the identity, and hence $p$ with section $e$ and multiplication $m$ is a fibrewise $H$-space.

\end{proof}

When dealing with fibrations, we have two different ways of generalizing rationalizations. The first is to work with the rationalization of the base and total space.
That is, given a fibration $p:Y\rightarrow B$ we can construct a rationalization $p_\Q:Y_\Q\rightarrow B_\Q$. This satisfies that for any commutative square

% https://q.uiver.app/#q=WzAsNCxbMCwwLCJZIl0sWzAsMSwiQiJdLFsxLDAsIlknIl0sWzEsMSwiQiciXSxbMCwxLCIiLDAseyJzdHlsZSI6eyJoZWFkIjp7Im5hbWUiOiJlcGkifX19XSxbMCwyXSxbMiwzLCIiLDIseyJzdHlsZSI6eyJoZWFkIjp7Im5hbWUiOiJlcGkifX19XSxbMSwzXV0=
\[\begin{tikzcd}[cramped]
	Y & {Y'} \\
	B & {B'}
	\arrow[from=1-1, to=1-2]
	\arrow[two heads, from=1-1, to=2-1]
	\arrow[two heads, from=1-2, to=2-2]
	\arrow[from=2-1, to=2-2]
\end{tikzcd}\]

where $Y',B'$ are rational spaces, there is a (unique up to homotopy) factorization through the rationalization. If we take this fibration and pull
back along the rationalization $p:B\rightarrow B_\Q$ we get a fibration which we will
denote $p^\Q:Y^\Q\twoheadrightarrow B$ which is called the `fibrewise rationalization' of $p$. It satisfies a universal property in the category of fibrations over $B$: given a map $f:Y\rightarrow R$ over $B$, where $R\twoheadrightarrow B$ is a fibration with fiber a rational space, $f$ factors uniquely up to homotopy through the fibrewise rationalization. In particular, the homotopy groups of the fiber are rationalizations of the homotopy groups of the fiber of $p$.

We will need the following lemma about fibrewise rationalizations.

\begin{thm}\label{ratfibheapstruc}
    Let $p:Y\twoheadrightarrow B$ be a fibration of simply connected spaces with simply connected fiber, with a relative minimal model $(M_B\otimes M_F, d)$ which is linear through dimension $k$. Then for any $n\leq k$ the $n^{th}$ Moore-Postnikov stage $L_n^\Q$ of the fibrewise rationalization $p^{\Q}:Y^{\Q}\twoheadrightarrow B$ is a fibrewise HM-space, and the maps $L_n^\Q\twoheadrightarrow L_{n-1}^\Q$ are all fibrewise HM-maps, as are the classifying maps $k_n:L_{n-1}\rightarrow B \times K(\pi_n(F)\otimes \Q, n+1)$.

\end{thm}

\begin{proof}
  We will proceed in two stages. We start by showing the rationalization $p_\Q:Y_\Q\twoheadrightarrow B_\Q$ is a fibrewise HM-space, and then we lift that structure to the fibrewise rationalization. To see that the rationalization is a fibrewise HM-space is straightforward: starting with the relative minimal model, we can construct a coMal'cev operation. In particular we have the relative minimal model presented as $M_B\otimes \wedge W_F$ where $W_F$ is the graded vector space built out of the rational homotopy groups of $F$. The linearity condition on the differential ensures that the map $M_B\otimes \wedge W_F\rightarrow (M_B\otimes \wedge W_F)^{\otimes_{\M_B}3}$ induced by the map sending $w\in W$ to $w\otimes 1\otimes 1-1\otimes w \otimes 1 + 1\otimes 1 \otimes w$ is a map of dgas. This map satisfies properties dual to those of the Mal'cev operation, and so it gives us a fibrewise heap structure on the rationalization. We want to first show that the k-invariants on each Postnikov stage of the rationalization are $HM$-maps. The spaces $K(\pi_n(F)\otimes \Q, n+1)$ are H-spaces, and so the operation $(a,b,c)=a-b+c$ endows it with the structure of an $HM$-space. This gives the trivial $B_\Q$ fibration $B_\Q \times K(\pi_n(F)\otimes \Q, n+1)$ the structure of a relative $HM$-space. A simple computation shows that the coMal'cev operation this induces on the relative minimal model on this fibration is precisely the map sending $\alpha$ to $\alpha\otimes 1 \otimes 1-1\otimes \alpha \otimes 1 + 1\otimes 1 \otimes \alpha$, from which it follows that the $k$-invariant on the corresponding stage of the Moore-Postnikov tower is an $HM$-map.

% https://q.uiver.app/#q=WzAsNCxbMCwwLCJMX3tuLFxcUX0iXSxbMCwxLCJMX3tuLTEsIFxcUX0iXSxbMSwwLCJCXFx0aW1lcyBFKFxccGlfbihGKVxcb3RpbWVzIFxcUSwgbikiXSxbMSwxLCJCXFx0aW1lcyBLKFxccGlfbihGKVxcb3RpbWVzIFxcUSwgbisxKSJdLFswLDEsIiIsMCx7InN0eWxlIjp7ImhlYWQiOnsibmFtZSI6ImVwaSJ9fX1dLFswLDJdLFsxLDMsImtfe24sXFxRfSJdLFsyLDNdLFswLDMsIiIsMSx7InN0eWxlIjp7Im5hbWUiOiJjb3JuZXIifX1dXQ==
\[\begin{tikzcd}[cramped]
	{L_{n,\Q}} & {B\times E(\pi_n(F)\otimes \Q, n)} \\
	{L_{n-1, \Q}} & {B\times K(\pi_n(F)\otimes \Q, n+1)}
	\arrow[two heads, from=1-1, to=2-1]
	\arrow[from=1-1, to=1-2]
	\arrow["{k_{n,\Q}}", from=2-1, to=2-2]
	\arrow[from=1-2, to=2-2]
	\arrow["\lrcorner"{anchor=center, pos=0.125}, draw=none, from=1-1, to=2-2]
  \end{tikzcd}\]

Because the maps on the bottom and right of this pullback square are $HM$-maps, so is the map $L_{n,\Q}\rightarrow L_{n-1, \Q}$. Them we need to show that this structure pulls back appropriately to the fibrewise rationalization. In particular since the fibrewise rationalization is a pullback, we can define the Mal'cev operation on $L_{n}^{\Q}$ by pulling back the operation on the rationalization:

% https://q.uiver.app/#q=WzAsNixbMiwyLCJCX1xcUSJdLFsxLDIsIkIiXSxbMSwxLCJMX25ee1xcUX0iXSxbMiwxLCJMX3tuLFxcUX0iXSxbMiwwLCJMX3tuLFxcUX1ee1xcdGltZXNfe0JfXFxRfSAzfSJdLFswLDAsIihMX3tufV57XFxRfSlee1xcdGltZXNfQiAzfSJdLFsxLDBdLFsyLDNdLFsyLDFdLFszLDBdLFsyLDAsIiIsMSx7InN0eWxlIjp7Im5hbWUiOiJjb3JuZXIifX1dLFs0LDMsIlxcdGF1X3tuLFxcUX0iXSxbNSwxLCIiLDAseyJjdXJ2ZSI6Mn1dLFs1LDRdLFs1LDIsIlxcdGF1X24iLDAseyJzdHlsZSI6eyJib2R5Ijp7Im5hbWUiOiJkYXNoZWQifX19XV0=
\[\begin{tikzcd}[cramped]
	{(L_{n}^{\Q})^{\times_B 3}} && {L_{n,\Q}^{\times_{B_\Q} 3}} \\
	& {L_n^{\Q}} & {L_{n,\Q}} \\
	& B & {B_\Q}
	\arrow[from=3-2, to=3-3]
	\arrow[from=2-2, to=2-3]
	\arrow[from=2-2, to=3-2]
	\arrow[from=2-3, to=3-3]
	\arrow["\lrcorner"{anchor=center, pos=0.125}, draw=none, from=2-2, to=3-3]
	\arrow["{\tau_{n,\Q}}", from=1-3, to=2-3]
	\arrow[curve={height=12pt}, from=1-1, to=3-2]
	\arrow[from=1-1, to=1-3]
	\arrow["{\tau_n}", dashed, from=1-1, to=2-2]
\end{tikzcd}\]

The associativity and Mal'cev properties are satisfied trivially since this is a pullback of a fibrewise $HM$-space structure on the rationalization.

\end{proof}

Finally, we will need one more technical lemma about fibrewise $H$-spaces. This is a fibrewise version of lemma $3.6$ in \cite{manin2022rational} which tells us that the multiplication by $r$ map on a fibrewise $H$ space kills torsion elements of cohomology.

\begin{lm}\label{multiplytokilltorsion}
    Let $p:H\rightarrow B$ be a fibrewise H-space of finite type, $A$ a finitely generated coefficient group, and $\alpha\in H^{(n)}(H; A)$ a cohomology class with the property that $t\alpha\in p^{*}(H^{n}(B; A))$ for some positive $t$. Then there is an $r> 0$ such that $\chi_r^*\alpha\in p^{*}(H^{n}(B; A))$, where $\chi_r:H\rightarrow H$ is the `multiplication by $r$ map,' i.e. $\chi_2(a)=m(a,a)$, and $\chi_n(a)=m(a, \chi_{n-1}(a))$.
\end{lm}

\begin{proof}

  %We first note that since the inclusion of the fiber $i:F\hookrightarrow H$ composes with the fibration to the trivial map, the pullback $i^{*}\alpha$ is a torsion class in $H^{n}(F,A)$.

  Note that in the language of the Serre spectral sequence we can rewrite this condition as saying that we have an element in $\alpha\in \bigoplus_{p+q=n}E^{p,q}_{\infty}$ such that $t\alpha\in E_{\infty}^{n,0}$ for some $t$. Then it suffices to show the following: suppose we have some torsion element of $\beta\in E^{p,q}_{\infty}$ where $q\geq 1$, then there exists some $r$ such that $\chi_{r}^{*}\alpha=0$. Indeed, since the direct sum above is finite, then we can take the direct sum decomposition of the element $\alpha$ into a finite sum of terms $\beta_{q}$ in $E_{p,q}^{\infty}$, for which some $r_{q}$ will suffice, and then $\chi_{r_{1}\dots r_{n}}$ will kill all but $\beta_{0}$ as desired.
  Recall that the $E_{1}$ page of the Serre spectral sequence has terms $E_{1}^{p,q}=C^{p}(B; H^{q}(F; A))$.

  Suppose then that we have some torsion cohomology class $\beta$ as above, there is some corresponding cocycle $\gamma$ which survives to $E^{p,q}_{\infty}$, but $t\gamma$ does not survive. By lemma $3.5$ of  \cite{manin2022rational} we know that there is some $s$ such that $\chi_{s}^{*}(H^{q}(F;A))\subseteq tH^{q}(F;A)$, from which we can conclude that $\chi_{s}^{*}(C^{p}(B; H^{q}(F; A)))\subseteq tC^{p}(B; H^{q}(F; A))$, and hence $\chi_{s}^{*}(\gamma)$ does not survive to $E^{p,q}_{\infty}$, as desired.

\end{proof}

\section{The Lifting Algorithm}\label{liftingalgsection}

We turn now to the lifting algorithm. In particular, we have the following theorem:
\begin{thm}\label{liftingalg}
    There is an algorithm that on input a diagram:
    % https://q.uiver.app/?q=WzAsMyxbMCwxLCJYIl0sWzEsMSwiQiJdLFsxLDAsIlkiXSxbMiwxLCJwIiwwLHsic3R5bGUiOnsiaGVhZCI6eyJuYW1lIjoiZXBpIn19fV0sWzAsMSwiZiIsMl1d

    \[\begin{tikzcd}
	& Y \\
	X & B
	\arrow["p", two heads, from=1-2, to=2-2]
	\arrow["f"', from=2-1, to=2-2]
\end{tikzcd}\]
and a relative minimal model $(M_B\otimes M_F, d_L)$ for $p:Y\twoheadrightarrow p$ where each of the spaces, and the fiber of $p$ is a simply connected finite type CW complex, and the minimal model is linear through the dimension of $X$, decides whether there exists a lift $g:X\rightarrow Y$ of $f$.

\end{thm}

\begin{proof}
    Let $d$ be the dimension of $X$. To simplify the proof a bit, we pullback $p$ across $f$ to obtain a fibration $\hat{p}:\hat{Y}\rightarrow X$, for which we will determine if a section exists.
    
    We will denote by $P_n$ the $n^{th}$ Moore-Postnikov stage of $\hat{p}:\hat{Y}\twoheadrightarrow X$. The linear relative minimal model pulls back to model $\hat{p}$ so we still have a relative minimal model for this fibration which is linear through the dimension of $X$. Using this minimal model, we can apply lemma \ref{dgalift} to construct a map of dgas dual to a lift. Here the algorithm might fail to find a lift, in which case we know none exists. Indeed, if a lift existed, applying the equivalence would produce a dual map on dgas, so that if no such map exists, no section can exist.

    We assume then that we found such a map. In particular we have $\epsilon:\mathcal{M}_X\otimes \wedge W_F\rightarrow \mathcal{M}_X$ which is a map of dgas, where $\mathcal{M}_X\otimes \wedge W_F$ with differential $d_L^*$ is the relative minimal model of $\hat{p}$.

    We will attempt to construct a section for $\hat{p}$ inductively as follows:
    for each $n$ through dimension $d$ we will construct a fibration $h_n:L_n\twoheadrightarrow X$, 
    as well as a section $e_n$ and rational equivalences $\phi_n:L_n\rightarrow P_n$, $\theta_{n}:P_{n}\rightarrow L_{n}$ over $X$, where the composition $\theta_{n}\circ \phi_{n}$ is multiplication by some integer $s_{n}$ under the fibrewise $H$-space structure on $L_{n}$ with the section $e_{n}$.
    In particular, these will be constructed so that each $h_n$ is a fibrewise $HM$-space with operation $\tau_n$, and there are maps $r_n:L_n\rightarrow L_{n-1}$ which commute with the fibrewise $HM$-structure and form a Moore-Postnikov tower for $h_d$.
    At each stage we will also fix a fibrewise rationalization $u_n$ to the corresponding Moore-Postnikov stage $L_n^\Q$ which is isomorphic to $P_n^\Q$. By lemma \ref{ratfibheapstruc} each $L_n^\Q$ has the structure of a fibrewise $HM$-space, and we will ensure that the maps $u_n$ are each $HM$-maps.

    At each stage, in the construction of $\phi_{n}$, it is possible for the construction to fail, and this indicates that no section of $\hat{p}$ exists. We will see why when we construct $\phi_{n}$.

    Since the space is simply connected, we set $h_1:X\twoheadrightarrow X$ to the identity and the other data is trivial.

    Then we assume we have constructed $h_{n-1}$, $e_{n-1}$, $\phi_{n-1}$,$\theta_{n-1}$ and $u_{n-1}$. We start by constructing the map $h_n$. We do this by picking $L_n\twoheadrightarrow L_{n-1}$ as a $K(\pi_n(F), n)$ fibration. This will then serve as the top stage on the Moore-Postnikov tower for $h_n:L_n\twoheadrightarrow X$. This is equivalent to picking a classifying map in $[L_{n-1}, K(\pi_n(F), n+1)]_B$, or alternatively, a cohomology class in $H^{n+1}(L_{n-1}; \pi_n(F))$.
    By the universal coefficient theorem this is isomorphic to
    \[
        \Hom(H_{n+1}(L_{n-1}), \pi_n(F))\oplus \Ext^1_\Z(H_n(L_{n-1}), \pi_n(F))
    \]

    From the relative minimal model of $\hat{p}$ we have a map 
    \[ 
        d_L^*|_{W_F^{(n)}}: W_F^{(n)}\rightarrow H^{n+1}(L_{n-1}; \Q)
    \]    
    and since $W_F^{(n)}$ is simply $(\pi_n(F)\otimes \Q)^*$ or equivalently $\Hom(\pi_n(F), \Q)$ we have an element of
    \[ 
        \Hom(\Hom(\pi_n(F), \Q), H^{(n+1)}(L_{n-1}; \Q))
    \]
    which we will denote $k_{n,\Q}$.

    Now since $\Ext(H_n(X), \Q)$ is trivial, the universal coefficient theorem lets us view the above as isomorphic to \[
    \Hom(\Hom(\pi_n(F), \Q), \Hom(H_{(n+1)}(L_{n-1}), \Q))
    \]
    by taking the dual map then, we obtain an element
    \[
    \Hom(H_{(n+1)}(L_{n-1})\otimes \Q, \pi_n(F)\otimes \Q)
    \]
    From this element $d_L^*|_{W_F^{(n)}}\in \Hom(H_{n+1}(L_{n-1})\otimes \Q, \pi_n(F)\otimes\Q)$ we will construct an element of $\Hom(H_{n+1}(L_{n-1}), \pi_n(F))$. In particular, there is a natural map $H_{n+1}(L_{n-1})\rightarrow H_{n+1}(L_{n-1}\otimes \Q)$ given by $a\mapsto a\otimes 1$, and so by composing with this map we can construct $\widetilde{d_L^*|_{W_F^{(n)}}}\in \Hom(H_{n+1}(L_{n-1}), \pi_n(F)\otimes \Q)$.
    Next we fix a minimal generating set for both $H_{n+1}(L_{n-1})$ and $\pi_n(F)$. In particular, from the basis for $W^{(n)}$ we have a minimal set of generators for the free part of $\pi_n(F)$ so we simply adjoin a minimal generating set for the torsion part.
    Next we consider $\widetilde{d_L^*|_{W_F^{(n)}}}(a)$ for each $a$ in the generating set for $H_{n+1}(L_{n-1})$, and we can write each of these as linear combinations of pure tensors $b\otimes\frac{p}{q}$ for $b$ in the minimal generating set of $\pi_n(F)$.
    Since the torsion part is killed in the tensor with $\Q$ we know that only elements of the free part of $\pi_n(F)$ show up in these linear combinations, and these we can identify with basis elements for $W^{(n)}$. We consider the collection of all $\frac{p}{q}$ arising as coefficients in these terms, and we can pick the least common multiple $Q$ of the $q$s.
    Then the subgroup $Q H_{n+1}(L_{n-1})$ lands in the image of $\pi_n(F)$ under the map $\pi_n(F)\rightarrow \pi_n(F)\otimes \Q$. Then precomposing with the multiplication by $Q\times N$ map gives us a map which lifts to $\tilde{k_n}:H_{n+1}(L_{n-1})\rightarrow \pi_n(F)$, where $N$ is an integer we will determine shortly.

    Since $H^{n+1}(L_{n-1};\pi_n(F))$ can be decomposed as $\Hom(H_{n+1}(L_{n-1}), \pi_n(F))\oplus Ext^1_\Z(H_n(L_{n-1}), \pi_n(F))$, we can simply perform the inclusion $\Hom(H_{n+1}(L_{n-1}), \pi_n(F))\hookrightarrow \Hom(H_{n+1}(L_{n-1}), \pi_n(F)) \oplus \Ext^1_\Z(H_n(L_{n-1}), \pi_n(F))$ on $\tilde{k_n}$ to get an element $k_n\in H^{n+1}(L_{n-1};\pi_n(F))$.

We now have a fibration $h_n:L_n\twoheadrightarrow X$. We will now construct the map $u_n$. This is simply fixing a rationalization along the dashed line in the square below:

% https://q.uiver.app/#q=WzAsNCxbMCwwLCJMX24iXSxbMSwwLCJMX25eXFxRIl0sWzAsMSwiTF97bi0xfSJdLFsxLDEsIkxfe24tMX1eXFxRIl0sWzAsMSwidV9uIiwwLHsic3R5bGUiOnsiYm9keSI6eyJuYW1lIjoiZGFzaGVkIn19fV0sWzAsMiwicl9uIiwyXSxbMSwzLCJyX25eXFxRIl0sWzIsMywidV97bi0xfSIsMl1d
\[\begin{tikzcd}
	{L_n} & {L_n^\Q} \\
	{L_{n-1}} & {L_{n-1}^\Q}
	\arrow["{u_n}", dashed, from=1-1, to=1-2]
	\arrow["{r_n}"', from=1-1, to=2-1]
	\arrow["{r_n^\Q}", from=1-2, to=2-2]
	\arrow["{u_{n-1}}"', from=2-1, to=2-2]
\end{tikzcd}\]

Next we construct the relative Mal'cev operation $\tau_n$.
To do this, we will consider the following diagram.

% https://q.uiver.app/#q=WzAsMTEsWzAsMCwiTF9uXFx0aW1lc19YIExfblxcdGltZXNfWCBMX24iXSxbMCwyLCJMX3tuLTF9XFx0aW1lc19YIExfe24tMX1cXHRpbWVzX1ggTF97bi0xfSJdLFsxLDEsIkxfbiJdLFsyLDAsIkxfbl5cXFEiXSxbNCwwLCJYXFx0aW1lcyBFKFxccGlfbihGKVxcb3RpbWVzXFxRLG4pIl0sWzQsMiwiWFxcdGltZXMgSyhcXHBpX24gKEYpXFxvdGltZXMgXFxRLCBuKzEpIl0sWzIsMiwiTF97bi0xfV5cXFEiXSxbMSwzLCJMX3tuLTF9Il0sWzMsMywiWFxcdGltZXMgSyhcXHBpX24oRiksbisxKSJdLFszLDEsIlhcXHRpbWVzIEUoXFxwaV9uKEYpLG4pIl0sWzMsMl0sWzAsMSwiaF9uXntcXHRpbWVzX1ggM30iLDJdLFswLDMsIlxcdGF1X25eXFxRXFxjaXJjIHVfbl57XFx0aW1lc19YIDN9Il0sWzMsNCwiXFxoYXR7a31fbl5cXFEiXSxbNCw1XSxbNiw1LCJrX25eXFxRIiwyXSxbMyw2XSxbMiw3LCJoX24iLDJdLFs3LDgsImtfbiIsMl0sWzEsNywiXFx0YXVfe24tMX0iLDJdLFs5LDhdLFs4LDVdLFs5LDRdLFsyLDMsInVfbiJdLFs3LDYsInVfe24tMX0iXSxbMiw5LCJcXGhhdHtrfV9uIl0sWzIsOCwiIiwxLHsic3R5bGUiOnsibmFtZSI6ImNvcm5lciJ9fV0sWzAsMiwiXFx0YXVfbiIsMSx7InN0eWxlIjp7ImJvZHkiOnsibmFtZSI6ImRhc2hlZCJ9fX1dLFszLDEwLCIiLDIseyJzdHlsZSI6eyJuYW1lIjoiY29ybmVyIn19XV0=
\[\begin{tikzcd}[cramped]
	{L_n\times_X L_n\times_X L_n} && {L_n^\Q} && {X\times E(\pi_n(F)\otimes\Q,n)} \\
	& {L_n} && {X\times E(\pi_n(F),n)} \\
	{L_{n-1}\times_X L_{n-1}\times_X L_{n-1}} && {L_{n-1}^\Q} & {} & {X\times K(\pi_n (F)\otimes \Q, n+1)} \\
	& {L_{n-1}} && {X\times K(\pi_n(F),n+1)}
	\arrow["{\tau_n^\Q\circ u_n^{\times_X 3}}", from=1-1, to=1-3]
	\arrow["{\tau_n}"{description}, dashed, from=1-1, to=2-2]
	\arrow["{h_n^{\times_X 3}}"', from=1-1, to=3-1]
	\arrow["{\hat{k}_n^\Q}", from=1-3, to=1-5]
	\arrow[from=1-3, to=3-3]
	\arrow["\lrcorner"{anchor=center, pos=0.125}, draw=none, from=1-3, to=3-4]
	\arrow[from=1-5, to=3-5]
	\arrow["{u_n}", from=2-2, to=1-3]
	\arrow["{\hat{k}_n}", from=2-2, to=2-4]
	\arrow["{h_n}"', from=2-2, to=4-2]
	\arrow["\lrcorner"{anchor=center, pos=0.125}, draw=none, from=2-2, to=4-4]
	\arrow[from=2-4, to=1-5]
	\arrow[from=2-4, to=4-4]
	\arrow["{\tau_{n-1}}"', from=3-1, to=4-2]
	\arrow["{k_n^\Q}"', from=3-3, to=3-5]
	\arrow["{u_{n-1}}", from=4-2, to=3-3]
	\arrow["{k_n}"', from=4-2, to=4-4]
	\arrow[from=4-4, to=3-5]
  \end{tikzcd}\]

Our aim is to construct a $\tau_n$ along the dashed line making the diagram commute. Essentially, we want to lift the Mal'cev operation from the fibrewise rationalization. We note that if we construct a map $\tau_n$ making the above diagram commute, it will commute with the map $u_n$.

$L_n$ is a pullback, so we can build $\tau_n$ by constructing a map to $X\times E(\pi_n(F), n)$ that commutes with the rest of the pullback square. Starting at $L_n\times_X L_n\times_X L_n$ we can follow along the map $\hat{k}_n^\Q\circ \tau_n^\Q\circ u_n^{\times_X 3}$. We want to lift this to a map to $X\times E(\pi_n(F), n+1)$, in such a way that it is also a lift of the map $k_n\circ \tau_{n-1}\circ h_n^{\times_X 3}:L_n^{\times_X 3}\rightarrow X\times K(\pi_n(F), n+1)$.
 We note that $X\times E(\pi_n(F), n)$ is also an $HM$-space. Then we have a map $\hat{\tau}_n:X\times E(\pi_n(F), n)^{3}\rightarrow X\times E(\pi_n(F), n)$, and composing this with $\hat{k}_n^{\times_X 3}$ provides a map $\psi_n:L_n^{\times_X 3} \rightarrow X \times E(\pi_n(F), n)$. Since $\hat{k}_n^\Q$ commutes with the fibrewise $HM$-space structure, we know that $\psi_n$ commutes with the diagram:

% https://q.uiver.app/#q=WzAsOSxbMCwwLCJMX25cXHRpbWVzX1ggTF9uXFx0aW1lc19YIExfbiJdLFswLDIsIkxfe24tMX1cXHRpbWVzX1ggTF97bi0xfVxcdGltZXNfWCBMX3tuLTF9Il0sWzEsMSwiTF9uIl0sWzIsMCwiTF9uXlxcUSJdLFs0LDAsIlhcXHRpbWVzIEUoXFxwaV9uKEYpXFxvdGltZXNcXFEsbikiXSxbNCwyLCJYXFx0aW1lcyBLKFxccGlfbiAoRilcXG90aW1lcyBcXFEsIG4rMSkiXSxbMiwyLCJMX3tuLTF9XlxcUSJdLFsxLDMsIkxfe24tMX0iXSxbMywxLCJYXFx0aW1lcyBFKFxccGlfbihGKSxuKSJdLFswLDEsImhfbl57XFx0aW1lc19YIDN9IiwyXSxbMCwzLCJcXHRhdV9uXlxcUVxcY2lyYyB1X25ee1xcdGltZXNfWCAzfSJdLFszLDQsIlxcaGF0e2t9X25eXFxRIl0sWzQsNV0sWzYsNSwia19uXlxcUSIsMl0sWzMsNl0sWzMsNSwiIiwyLHsic3R5bGUiOnsibmFtZSI6ImNvcm5lciJ9fV0sWzIsNywiaF9uIiwyXSxbMSw3LCJcXHRhdV97bi0xfSIsMl0sWzgsNF0sWzIsMywidV9uIl0sWzcsNiwidV97bi0xfSJdLFsyLDgsIlxcaGF0e2t9X24iXSxbMCwyLCJcXHRhdV9uIiwxLHsic3R5bGUiOnsiYm9keSI6eyJuYW1lIjoiZGFzaGVkIn19fV0sWzAsOCwiXFxwc2lfbiIsMCx7ImN1cnZlIjotNX1dXQ==
\[\begin{tikzcd}
	{L_n\times_X L_n\times_X L_n} && {L_n^\Q} && {X\times E(\pi_n(F)\otimes\Q,n)} \\
	& {L_n} && {X\times E(\pi_n(F),n)} \\
	{L_{n-1}\times_X L_{n-1}\times_X L_{n-1}} && {L_{n-1}^\Q} && {X\times K(\pi_n (F)\otimes \Q, n+1)} \\
	& {L_{n-1}}
	\arrow["{h_n^{\times_X 3}}"', from=1-1, to=3-1]
	\arrow["{\tau_n^\Q\circ u_n^{\times_X 3}}", from=1-1, to=1-3]
	\arrow["{\hat{k}_n^\Q}", from=1-3, to=1-5]
	\arrow[from=1-5, to=3-5]
	\arrow["{k_n^\Q}"', from=3-3, to=3-5]
	\arrow[from=1-3, to=3-3]
	\arrow["\lrcorner"{anchor=center, pos=0.125}, draw=none, from=1-3, to=3-5]
	\arrow["{h_n}"', from=2-2, to=4-2]
	\arrow["{\tau_{n-1}}"', from=3-1, to=4-2]
	\arrow[from=2-4, to=1-5]
	\arrow["{u_n}", from=2-2, to=1-3]
	\arrow["{u_{n-1}}", from=4-2, to=3-3]
	\arrow["{\hat{k}_n}", from=2-2, to=2-4]
	\arrow["{\tau_n}"{description}, dashed, from=1-1, to=2-2]
	\arrow["{\psi_n}", curve={height=-50pt}, from=1-1, to=2-4]
\end{tikzcd}\]

% curve on last line edited manually to prevent intersections

and so it remains to show that it commutes with

% https://q.uiver.app/#q=WzAsNSxbMCwwLCJMX25cXHRpbWVzX1ggTF9uXFx0aW1lc19YIExfbiJdLFswLDIsIkxfe24tMX1cXHRpbWVzX1ggTF97bi0xfVxcdGltZXNfWCBMX3tuLTF9Il0sWzEsMywiTF97bi0xfSJdLFszLDEsIlhcXHRpbWVzIEUoXFxwaV9uKEYpLG4pIl0sWzMsMywiWFxcdGltZXMgSyhcXHBpX24oRiksbisxKSJdLFswLDEsImhfbl57XFx0aW1lc19YIDN9IiwyXSxbMSwyLCJcXHRhdV97bi0xfSIsMl0sWzIsNCwia19uIiwyXSxbMyw0XSxbMCwzLCJcXHBzaV9uIiwwLHsiY3VydmUiOi0zfV1d
\[\begin{tikzcd}
	{L_n\times_X L_n\times_X L_n} \\
	&&& {X\times E(\pi_n(F),n)} \\
	{L_{n-1}\times_X L_{n-1}\times_X L_{n-1}} \\
	& {L_{n-1}} && {X\times K(\pi_n(F),n+1)}
	\arrow["{h_n^{\times_X 3}}"', from=1-1, to=3-1]
	\arrow["{\tau_{n-1}}"', from=3-1, to=4-2]
	\arrow["{k_n}"', from=4-2, to=4-4]
	\arrow[from=2-4, to=4-4]
	\arrow["{\psi_n}", curve={height=-18pt}, from=1-1, to=2-4]
\end{tikzcd}\]

which is where we determine $N$ as mentioned above. Commutativity of the diagram above hinges on simply the two maps to $X\times K(\pi_n(F), n+1)$ agreeing, or equivalently equality of the pair of cohomology classes in $H^{(n+1)}(L_n^{\times_X 3}; \pi_n(F))$. One of these maps factors through $E(\pi_n(F), n)$ and so the corresponding cohomology class is trivial. Then we only need to look at the cohomology class from the bottom path in the diagram. Since the diagram commutes after rationalization, we know this class is a torsion element. If we construct a $\tilde{k}_n$ with $N=1$ then look at the corresponding torsion class, we set $N$ to be $m$ times the order of this class, where $m$ will be determined by an analogous argument for extending the section. Then considering the classes $\kappa_n$ and $\tilde{\kappa_n}$ in $H^{(n+1)}(L_{n-1}; \pi_n(F))$ represented by $k_n$ and $\tilde{\kappa_n}$, we have $\kappa_n = N\tilde{\kappa_n}$,
and so pulling back across the map $\tau_{n-1}\circ r_n^{\times_X 3}$
we get that the torsion element from $\tilde{\kappa}_n$ will be killed by this multiplication. Then the class we get pulling back $k_n$ is trivial, and so the above diagram commutes. Then we can pull back $\psi_n$ across $k_n$ and we have a $\tau_n$ as desired.

That $\tau_n$ satisfies the conditions for a fibrewise Mal'cev operation is a straightforward consequence of the fact that we are defining $\tau_n$ by pulling back a Mal'cev operation across a map that respects the fibrewise Mal'cev operation.

Now we construct the section.

Again we are looking for a map along the dashed line in the following diagram:
% https://q.uiver.app/#q=WzAsOSxbMCwyLCJYIl0sWzEsMSwiTF9uIl0sWzIsMCwiTF9uXlxcUSJdLFs0LDAsIlhcXHRpbWVzIEUoXFxwaV9uKEYpXFxvdGltZXNcXFEsbikiXSxbNCwyLCJYXFx0aW1lcyBLKFxccGlfbiAoRilcXG90aW1lcyBcXFEsIG4rMSkiXSxbMiwyLCJMX3tuLTF9XlxcUSJdLFsxLDMsIkxfe24tMX0iXSxbMywzLCJYXFx0aW1lcyBLKFxccGlfbihGKSxuKzEpIl0sWzMsMSwiWFxcdGltZXMgRShcXHBpX24oRiksbikiXSxbMiwzLCJcXGhhdHtrfV9uXlxcUSJdLFszLDRdLFs1LDQsImtfbl5cXFEiLDJdLFsyLDVdLFsyLDQsIiIsMix7InN0eWxlIjp7Im5hbWUiOiJjb3JuZXIifX1dLFsxLDYsImhfbiIsMl0sWzYsNywia19uIiwyXSxbMCw2LCJlX3tuLTF9IiwyXSxbOCw3XSxbNyw0XSxbOCwzXSxbMSwyLCJ1X24iXSxbNiw1LCJ1X3tuLTF9Il0sWzEsOCwiXFxoYXR7a31fbiJdLFsxLDcsIiIsMSx7InN0eWxlIjp7Im5hbWUiOiJjb3JuZXIifX1dLFswLDIsImVfbl5cXFEiLDIseyJjdXJ2ZSI6LTR9XSxbMCwxLCJlX24iLDIseyJzdHlsZSI6eyJib2R5Ijp7Im5hbWUiOiJkYXNoZWQifX19XV0=
\[\begin{tikzcd}[cramped]
	&& {L_n^\Q} && {X\times E(\pi_n(F)\otimes\Q,n)} \\
	& {L_n} && {X\times E(\pi_n(F),n)} \\
	X && {L_{n-1}^\Q} && {X\times K(\pi_n (F)\otimes \Q, n+1)} \\
	& {L_{n-1}} && {X\times K(\pi_n(F),n+1)}
	\arrow["{\hat{k}_n^\Q}", from=1-3, to=1-5]
	\arrow[from=1-5, to=3-5]
	\arrow["{k_n^\Q}"', from=3-3, to=3-5]
	\arrow[from=1-3, to=3-3]
	\arrow["\lrcorner"{anchor=center, pos=0.125, rotate=45}, draw=none, from=1-3, to=3-5]
	\arrow["{h_n}"', from=2-2, to=4-2]
	\arrow["{k_n}"', from=4-2, to=4-4]
	\arrow["{e_{n-1}}"', from=3-1, to=4-2]
	\arrow[from=2-4, to=4-4]
	\arrow[from=4-4, to=3-5]
	\arrow[from=2-4, to=1-5]
	\arrow["{u_n}", from=2-2, to=1-3]
	\arrow["{u_{n-1}}", from=4-2, to=3-3]
	\arrow["{\hat{k}_n}", from=2-2, to=2-4]
	\arrow["\lrcorner"{anchor=center, pos=0.125}, draw=none, from=2-2, to=4-4]
	\arrow["{e_n^\Q}"', curve={height=-24pt}, from=3-1, to=1-3]
	\arrow["{e_n}"', dashed, from=3-1, to=2-2]
\end{tikzcd}\]

and by an analogous argument to the one above, the obstructions to making such a lift lie in the torsion part of $H^{n+1}(X; \pi_{n}(F))$, and so constructing this diagram with $m=1$ above will give us such a torsion element, and setting $m$ to be the order of this element will kill the obstruction, allowing us to pick an $e_{n}$.

We are now ready to construct $\phi_n:L_n\rightarrow P_n$. By construction, $L_{n}$ is fibrewise rationally equivalent to $P_{n}$, since the fibers are rationally equivalent, and the $k$-invariants are the same up to torsion. Then it remains only to show that we can actually compute such a rational equivalence. Suppose we try to build a map along the dashed line in the following diagram, making it commute:
% https://q.uiver.app/#q=WzAsNixbMCwwLCJMX24iXSxbMCwxLCJMX3tuLTF9Il0sWzEsMSwiUF97bi0xfSJdLFsxLDAsIlBfbiJdLFsyLDAsIkUoXFxwaV9uKEYpLCBuKSJdLFsyLDEsIksoXFxwaV9uKEYpLCBuKzEpIl0sWzAsMSwicl9uIiwyXSxbMSwyLCJcXHBoaV97bi0xfSJdLFszLDJdLFszLDRdLFsyLDVdLFs0LDVdLFszLDUsIiIsMSx7InN0eWxlIjp7Im5hbWUiOiJjb3JuZXIifX1dLFswLDMsIlxccGhpX24iLDAseyJzdHlsZSI6eyJib2R5Ijp7Im5hbWUiOiJkYXNoZWQifX19XV0=
\[\begin{tikzcd}[cramped]
	{L_n} & {P_n} & {E(\pi_n(F), n)} \\
	{L_{n-1}} & {P_{n-1}} & {K(\pi_n(F), n+1)}
	\arrow["{r_n}"', from=1-1, to=2-1]
	\arrow["{\phi_{n-1}}", from=2-1, to=2-2]
	\arrow[from=1-2, to=2-2]
	\arrow[from=1-2, to=1-3]
	\arrow[from=2-2, to=2-3]
	\arrow[from=1-3, to=2-3]
	\arrow["\lrcorner"{anchor=center, pos=0.125}, draw=none, from=1-2, to=2-3]
	\arrow["{\phi_n}", dashed, from=1-1, to=1-2]
  \end{tikzcd}\]
Since the map $\phi_{n-1}$ is a rational equivalence, and $r_{n}$ followed by $k_{n}$ is $0$, we know that the obstruction to lifting $\phi_{n-1}\circ r_{n}$ is torsion. With the choice of section, we can endow $L_{n}$ with a fibrewise $H$-space structure, and then an application of lemma \ref{multiplytokilltorsion} tells us that if we precompose $r_{n}$ in the above diagram with $\chi_{k}$ for some integer $k$, we can push the obstruction to lifting to an element of $H^{n+1}(B; \pi_{n}(F))$, which will in particular be an obstruction to lifting $\phi_{n-1}\circ r_{n}\circ e_{n}$ to $P_{n}$. We claim here that either this obstruction vanishes, or there is no section to $P_{n}$.

Indeed, suppose there is a section $\gamma:X\rightarrow P_{n}$. Then consider the following diagram.

% https://q.uiver.app/#q=WzAsMTAsWzAsMiwiTF97bi0xfSJdLFsxLDIsIlBfe24tMX0iXSxbMiwyLCJMX3tuLTF9Il0sWzAsMywiSyhcXHBpX24oRiksbisxKSJdLFsxLDMsIksoXFxwaV9uKEYpLG4rMSkiXSxbMiwzLCJLKFxccGlfbihGKSxuKzEpIl0sWzEsMSwiUF9uIl0sWzIsMSwiTF9uIl0sWzAsMSwiTF9uIl0sWzEsMCwiWCJdLFswLDEsIlxccGhpX3tuLTF9Il0sWzEsMiwiXFx0aGV0YV97bi0xfSJdLFswLDMsImtfbiIsMl0sWzEsNCwia15wX24iLDJdLFsyLDUsImtfbiJdLFszLDQsIlxcYWxwaGEiLDJdLFs0LDUsIlxcYmV0YSIsMl0sWzYsMSwicF9uIiwyXSxbNywyLCJyX24iLDJdLFs4LDAsInJfbiIsMl0sWzksNiwiXFxnYW1tYSIsMl0sWzksNywiXFx0aWxkZXtcXGdhbW1hfSIsMSx7InN0eWxlIjp7ImJvZHkiOnsibmFtZSI6ImRhc2hlZCJ9fX1dXQ==
\[\begin{tikzcd}
	& X \\
	{L_n} & {P_n} & {L_n} \\
	{L_{n-1}} & {P_{n-1}} & {L_{n-1}} \\
	{K(\pi_n(F),n+1)} & {K(\pi_n(F),n+1)} & {K(\pi_n(F),n+1)}
	\arrow["\gamma"', from=1-2, to=2-2]
	\arrow["{\tilde{\gamma}}"{description}, dashed, from=1-2, to=2-3]
	\arrow["{r_n}"', from=2-1, to=3-1]
	\arrow["{p_n}"', from=2-2, to=3-2]
	\arrow["{r_n}"', from=2-3, to=3-3]
	\arrow["{\phi_{n-1}}", from=3-1, to=3-2]
	\arrow["{k_n}"', from=3-1, to=4-1]
	\arrow["{\theta_{n-1}}", from=3-2, to=3-3]
	\arrow["{k^p_n}"', from=3-2, to=4-2]
	\arrow["{k_n}", from=3-3, to=4-3]
	\arrow["\alpha"', from=4-1, to=4-2]
	\arrow["\beta"', from=4-2, to=4-3]
\end{tikzcd}\]

We want to show that the map $\gamma$ lifts to a section $\tilde{\gamma}$. First, we note that we have yet to define the maps $\alpha$ and $\beta$. These can be described by maps on $\pi_{n}(F)$, and in order to make the diagram commute, we simply have to make sure we pull back the correct cohomology classes on $L_{n}$ and $P_{n}$ in each spot. In particular, by construction of $k_{n}$ we can set $\alpha$ to be the map which kills torsion elements and multiplies the free part by $QN$. Similarly, $\beta$ can be the map which kills torsion and multiplies the free part by $\frac{s_{n}}{QN}$. Since $r_{n}$ is a pullback of the universal $K(\pi_{n}(F), n)$ fibration, whether $\theta_{n-1}\circ p_{n} \circ \gamma$ lifts to $L_{n}$ depends only on the map to $K(\pi_{n}(F), n+1)$, but by commutativity, this is simply $\beta\circ k_{n}^{P}\circ p_{n}\gamma$ and since this factors through the homotopically trivial map $k_{n}^{P}\circ p_{n}$ we know it is homotopically trivial. Then we have a section $\tilde{\gamma}:X \rightarrow L_{n}$. An identical argument on the other side of the diagram allows us to lift $\phi_{n-1}\circ r_{n}\circ \tilde{\gamma}$ to $P_{n}$, and create a section $\gamma ':X\rightarrow P_{n}$. Then we have the following commutative diagram:

% https://q.uiver.app/#q=WzAsNSxbMSwwLCJYIl0sWzAsMSwiTF9uIl0sWzIsMSwiUF9uIl0sWzAsMiwiTF97bi0xfSJdLFsyLDIsIlBfe24tMX0iXSxbMyw0LCJcXHBoaV97bi0xfSJdLFsxLDMsInJfbiJdLFsyLDQsInBfbiIsMl0sWzAsMSwiXFx0aWxkZXtcXGdhbW1hfSJdLFswLDIsIlxcZ2FtbWEnIiwyXSxbMSwyLCJcXHBoaV9uIiwxLHsic3R5bGUiOnsiYm9keSI6eyJuYW1lIjoiZGFzaGVkIn19fV1d
\[\begin{tikzcd}
	& X \\
	{L_n} && {P_n} \\
	{L_{n-1}} && {P_{n-1}}
	\arrow["{\tilde{\gamma}}", from=1-2, to=2-1]
	\arrow["{\gamma'}"', from=1-2, to=2-3]
	\arrow["{\phi_n}"{description}, dashed, from=2-1, to=2-3]
	\arrow["{r_n}", from=2-1, to=3-1]
	\arrow["{p_n}"', from=2-3, to=3-3]
	\arrow["{\phi_{n-1}}", from=3-1, to=3-3]
\end{tikzcd}\]

Using the same argument above, this time endowing $L_{n}$ with a fibrewise $H$-space structure with the section $\tilde{\gamma}$, lemma \ref{multiplytokilltorsion} together with the existence of $\gamma'$ guarantees that the obstruction to lifting $r_{n}\circ \phi_{n-1}$ to $\phi_{n}$ vanishes. This contradicts that $e_{n-1}\circ r_{n}\circ e_{n}$ does not lift since $e_{n}\circ \phi_{n}$ now provides a lift.

Finally then it remains to construct $\theta_{n}$. In particular we want to construct a $\theta_{n}$ so that $\theta_{n}\circ \phi_{n}$ is $\chi_{s_{n}}$ on $L_{n}$ for some integer $s_{n}$. For each such integer then we have a diagram:

% https://q.uiver.app/#q=WzAsOCxbMSwxLCJQX3tuLTF9Il0sWzIsMSwiTF97bi0xfSJdLFswLDAsIkxfbiJdLFswLDEsIlBfbiJdLFszLDEsIkxfe24tMX0iXSxbMywwLCJMX24iXSxbMiwyLCJLKFxccGlfbihGKSwgbisxKSJdLFszLDIsIksoXFxwaV9uKEYpLCBuKzEpIl0sWzIsMywiXFxwaGlfbiJdLFszLDAsInBfbiIsMl0sWzAsMSwiXFx0aGV0YV97bi0xfSIsMl0sWzEsNCwiXFxjaGlfe1N9IiwyXSxbNSw0XSxbMiw1LCJcXGNoaV97c19ufSJdLFszLDUsIlxcdGhldGFfbiIsMCx7InN0eWxlIjp7ImJvZHkiOnsibmFtZSI6ImRhc2hlZCJ9fX1dLFsxLDYsImtfbiIsMl0sWzQsN10sWzYsNywiXFxjaGlfUyIsMl1d
\[\begin{tikzcd}
	{L_n} &&& {L_n} \\
	{P_n} & {P_{n-1}} & {L_{n-1}} & {L_{n-1}} \\
	&& {K(\pi_n(F), n+1)} & {K(\pi_n(F), n+1)}
	\arrow["{\chi_{s_n}}", from=1-1, to=1-4]
	\arrow["{\phi_n}", from=1-1, to=2-1]
	\arrow[from=1-4, to=2-4]
	\arrow["{\theta_n}", dashed, from=2-1, to=1-4]
	\arrow["{p_n}"', from=2-1, to=2-2]
	\arrow["{\theta_{n-1}}"', from=2-2, to=2-3]
	\arrow["{\chi_{S}}"', from=2-3, to=2-4]
	\arrow["{k_n}"', from=2-3, to=3-3]
	\arrow["{k_{n}}", from=2-4, to=3-4]
	\arrow["{\chi_S}"', from=3-3, to=3-4]
\end{tikzcd}\]

In particular, the obstruction to lifting $\theta_{n-1}\circ p_{n}$ is given by the cohomology class $k_{n}\circ \theta_{n-1}\circ p_{n}$ and since $\theta_{n-1}$ is a rational equivalence, the obstruction to such a lift is torsion. Then for some $S$ the obstruction vanishes, and we can produce a lift $\theta_{n}$.

Having completed our induction then, we can simply compose the final section $e_{n}$ and rational equivalence $\phi_{n}$ to obtain a section for $\hat{p}$ as desired.

\end{proof}

\section{Proof of Theorem \ref{mainresult}}\label{mainresultsection}

We are now ready to put together the proof of the main result. Suppose then we are given the triple $(M, N, f)$ as in the statement of the theorem. Then applying theorem \ref{universalliftsaremonos} we want to prove that the existence of such a lift is decidable.

There are two obstructions to using theorem \ref{liftingalg} then, firstly we need to construct a relative minimal model of the bundle $\Mono(m-\planes,n-\planes)\twoheadrightarrow BSO(m)\times BSO(n)$ which has linear differential and secondly need to address the possibility that $M$ is not simply connected, (preventing the use of any of the lifting algorithms.)

To construct the relative minimal model the first step is to determine a minimal model for both the base and the fiber. We consider two cases, based on the parity of $n$. In the case that $n$ is even then, the rational cohomology of $BSO(n)$ has a generator for each Pontrjagin class, and one for the Euler class which squares to the top Pontrjagin class. Since the codimension is odd, so is $m$ and hence $BSO(m)$ has rational cohomology generated only by the Pontrjagin classes. Then we have for the base the minimal model
\[\B=\Q\langle \alpha_i^{(4i)}, \beta_j^{(4j)}, \varepsilon^{(n)}\rangle \]
where $i\in\{1,...,\frac{n}{2}-1\}$, $j\in\{1,...,\frac{m-1}{2}\}$.

The fiber is the Stiefel manifold $V_m(\R^n)$ which is a homogeneous space $SO(n)/SO(n-m)$. To find the minimal model of this we will use the Cartan-Weil model for a homogeneous space, as in \cite{felix2008algebraic}. In particular this allows us to model $V_m(\R^n)$ up to homotopy via a fibration $SO(n)\hookrightarrow V_m(\R^n)\twoheadrightarrow BSO(n-m)$, which is given by the pullback of the universal $SO(n)$ fibration across the map $BSO(n-m)\rightarrow BSO(n)$ induced by the inclusion $SO(n-m)\hookrightarrow SO(n)$. Putting this together the underlying graded vector space generating the minimal model for the fiber is
\[V_F=\Q\{\gamma_k^{(4k-1)}, \sigma^{(n-1)}\}\]
with $k\in \{\frac{n-m+1}{2},...,\frac{n}{2}-1\}$. 

Then the relative minimal model $(\B\otimes\wedge V_F, \tilde d)$ is determined by the restriction of the differential to $V_F$. In order to compute this, we start by constructing a map 
% https://q.uiver.app/#q=WzAsMyxbMCwwLCJFU08obSlcXHRpbWVzIEVTTyhuKSJdLFsyLDAsIk1vbm8obS1wbGFuZXMsIG4tcGxhbmVzKSJdLFsxLDEsIkJTTyhtKVxcdGltZXMgQlNPKG4pIl0sWzAsMSwiZiJdLFswLDIsIiIsMix7InN0eWxlIjp7ImhlYWQiOnsibmFtZSI6ImVwaSJ9fX1dLFsxLDIsIiIsMCx7InN0eWxlIjp7ImhlYWQiOnsibmFtZSI6ImVwaSJ9fX1dXQ==
\[\begin{tikzcd}
	{ESO(m)\times ESO(n)} && {\Mono(m-\planes, n-\planes)} \\
	& {BSO(m)\times BSO(n)}
	\arrow["f", from=1-1, to=1-3]
	\arrow[two heads, from=1-1, to=2-2]
	\arrow[two heads, from=1-3, to=2-2]
\end{tikzcd}\]
which we will define as follows. Over a particular local trivialization $U_{i}$ of $ESO(m)\times ESO(n)\rightarrow BSO(m)\times BSO(n)$ we consider a point $(p,O_{m},O_{n})$ where $p\in U_{i}$, $O_{m}$ is an orthogonal $m\times m$ matrix and $O_{n}$ is an orthogonal $n\times n$ matrix. Then the columns of $O_{m}I_{m,n}O_{n}^{\dag}$ will provide an orthogonal frame, where $I_{m,n}$ is the $m\times n$ matrix sending the standard basis of $\R^{m}$ to the first $m$ elements of the standard basis of $\R^{n}$. To check that this defines a coherent map is simply a matter of checking that it behaves the same on an intersection of local trivializations, but this follows immediately from the definition of $\Mono(m-\planes, n-\planes)$.

% Figure out how explicit this argument needs to get.

In the world of $\Q$-dgas then this map is dual to a map over $\B$
\[(\B\otimes \wedge V_F, \tilde d)\overset{f^*}{\rightarrow} \B\langle a_i^{(4i-1)}, e^{(n-1)}, b_j^{(4j-1)}| da_i=\alpha_i, db_i=\beta_i, de=\varepsilon\rangle \]
and since this is a DGA map over $B$, it is determined by its action on elements of $V_F$. We then have to determine both $f^*(\gamma_k)$ and $f^*(\sigma)$. 

$\gamma_k$ is an element of $(\pi_{(4k-1)}(V_m(\R^n))\otimes \Q)^*$ and in this context $f^*$ is the dual of the map induced by $f$, $\pi_{(4k-1)}(SO(m)\times SO(n))\rightarrow \pi_{(4k-1)}(V_m(\R^n))$. To understand what $f$ does to the fiber $SO(m)\times SO(n)$ we decompose it as a sequence of steps: 
\[SO(m)\times SO(n)\overset {\text{include}\times {-1}}{\longrightarrow} SO(n)\times SO(n)\overset{\text{multiply}}{\longrightarrow} SO(n)\overset{\text{project}}{\longrightarrow} V_m(\R^n)\]

and dualizing this allows us to write

\[ f^*(\gamma_k)=\begin{cases*}
        a_k-b_k & $k\leq m$ \\
        a_k & otherwise 
    \end{cases*}
\]
and
\[ f^*(\sigma)=\varepsilon\]
which allows us to determine $f^*$. Because this map has to commute with the differential, we can conclude that for $k\leq m$
\[f^*(\tilde d(\gamma_k))=\beta_k-\alpha_k\]
but the only preimage of $\beta_k-\alpha_k$ under $f^*$ is $\beta_k-\alpha_k$ and so $d(\gamma_k)=\beta_k-\alpha_k$. Similarly for $k>m$ we conclude $\tilde d(\gamma_k)=a_k$ and $\tilde d(\sigma)=\varepsilon$. Then we have computed a relative minimal model $(\B\otimes \wedge V_F, \tilde d)$ for the bundle in our lifting problem.

In the case where $n$ is odd, and so $m$ is even, we have for the base
\[
  \mathcal{B}=\Q\langle \alpha^{(4i)}, \beta^{(4j)}, \varepsilon(m) \rangle
\]
with $i\in \{1,\dots, \frac{n-1}{2}\}, j\in \{1,\dots \frac{m}{2}-1 \}$. For the fiber we have the vector space \[V_{F}=\Q\{\gamma_{k}^{4k-1}\}\]
where $k\in \{\frac{n-m+1}{2},\dots, \frac{n-1}{2}\}$ (the argument is the same as the other case, we simply don't have any Euler class in the fiber since both $n$ and $n-m$ are odd.) Again we construct the map $f$ as before, and it still dualizes a map which decomposes on the fiber in the same way.

Finally, we address the case where $M$ is not simply connected. Essentially we are going to replace $M$ by a simply connected complex and create a lifting problem so that a lift exists exactly when one exists over $M$. Since the construction is nearly identical to the plus construction introduced in \cite{kervaire1969smooth}, we will call this space $M^{+}$ (the only difference is that for the plus construction we want a space with perfect fundamental group so we get a space with actually identical cohomology).

We construct $M^{+}$ in two steps, following fairly directly the idea for the plus construction. First, we pick a generating set for $\pi_{1}(M)$, (for instance choosing the complement of a spanning tree of the $1$-skeleton.) We then add $2$-cells with attaching maps along each such generator. We call this space $\tilde{M}$, and note that it is simply connected. We then consider the homology long exact sequence of the pair $(\tilde{M}, M)$. In particular since the space $\tilde{M}/M$ is a bouquet of $2$-spheres we have the sequence
\[
    0\rightarrow H_{2}(M)\rightarrow H_{2}(\tilde{M})\rightarrow H_{2}(\tilde{M}, M)\overset{\delta}{\rightarrow} H_{1}(M)\rightarrow 0
\]

Which gives the short exact sequence

\[
0\rightarrow H_{2}(M)\rightarrow H_{2}(\tilde{M})\rightarrow \ker\delta\rightarrow 0
\]

Since $H^{2}(\tilde{M}, M)$ is free abelian, so is $\ker\delta$. Then $H_{2}(\tilde{M})$ decomposes as a direct sum $H_{2}(M)\oplus F$ for a free abelian group $F$. Since $\tilde{M}$ is simply connected, the Hurewicz homomorphism gives us that each element of $H_{2}(M)$ is represented by a map $S^2\rightarrow \tilde{M}$. Then to obtain $M^{+}$ we attach $3$-cells with attaching maps representing a basis of $F$. Since the attaching maps form a basis of a subgroup of the free part of $H^{2}(\tilde{M})$ when we look at the homology sequence of the pair $(M^{+}, \tilde{M})$ the map $H_{3}(M^{+},\tilde{M})\rightarrow H_{2}(\tilde{M})$ has trivial kernel, and so the homology groups above degree $2$ are all isomorphic. In particular then we have a space $M^{+}$ which is simply connected but in degree $2$ and higher has isomorphic homology to $M$.

Then consider the following diagram:

% https://q.uiver.app/#q=WzAsNCxbMCwxLCJNIl0sWzAsMiwiTV4rIl0sWzEsMSwiQlNPKG0pXFx0aW1lcyBCU08obikiXSxbMSwwLCJNb25vKG0tcGxhbmVzLG4tcGxhbmVzKSJdLFswLDEsIiIsMCx7InN0eWxlIjp7InRhaWwiOnsibmFtZSI6Imhvb2siLCJzaWRlIjoidG9wIn19fV0sWzAsMl0sWzMsMiwiIiwwLHsic3R5bGUiOnsiaGVhZCI6eyJuYW1lIjoiZXBpIn19fV0sWzAsMywiZiIsMCx7InN0eWxlIjp7ImJvZHkiOnsibmFtZSI6ImRhc2hlZCJ9fX1dLFsxLDIsImciLDAseyJzdHlsZSI6eyJib2R5Ijp7Im5hbWUiOiJkYXNoZWQifX19XSxbMSwzLCJcXHRpbGRle2Z9IiwxLHsibGFiZWxfcG9zaXRpb24iOjcwLCJvZmZzZXQiOi0xLCJzdHlsZSI6eyJib2R5Ijp7Im5hbWUiOiJkYXNoZWQifX19XV0=
\[\begin{tikzcd}[cramped]
	& {Mono(m-planes,n-planes)} \\
	M & {BSO(m)\times BSO(n)} \\
	{M^+}
	\arrow[two heads, from=1-2, to=2-2]
	\arrow["f", dashed, from=2-1, to=1-2]
	\arrow[from=2-1, to=2-2]
	\arrow[hook, from=2-1, to=3-1]
	\arrow["{\tilde{f}}"{description, pos=0.7}, shift left, dashed, from=3-1, to=1-2]
	\arrow["g", dashed, from=3-1, to=2-2]
  \end{tikzcd}\]

We construct the lifting problem according to lemma \ref{universalliftsaremonos}, and we want to decide if an $f$ exists. Note that since $BSO(m)\times BSO(n)$ is rationally an $H$-space, we can use the main result of \cite{manin2022rational} to construct an extension. We can pick a $g$ and then look for a lift $\tilde{f}$ and if such a lift exists we are done since restricting to $M$ provides a lift $f$. Then suppose no such lift exists. Then in particular no such lift exists rationally, and the first obstruction to finding such a lift lies in $H^{n}(M^{+}; \pi_{n-1}(F)\otimes \Q)$ where $F$ is the fiber $V_{m}(\R^{n})$ of the fibration over $BSO(m)\times BSO(n)$. By the universal coefficient theorem, $H^{n}(M^{+}; \pi_{n-1}(F)\otimes \Q)\cong \Hom(H_{n}(M^{+}), \pi_{n-1}(F)\otimes \Q)$ and since the inclusion $i:M\hookrightarrow M^{+}$ induces an isomorphism between $H_{n}(M^{+})$ and $H_{n}(M)$ above degree $2$ (and $H_{1}(M^{+})$ is trivial) we can conclude that the obstruction is an obstruction to lifting $f$ as well.

To conclude then we summarize the steps of the algorithm:

{\bf Input}:
\begin{itemize}
  \item A pair of closed oriented smooth manifolds, $M$,$N$ as a pair of simplicial complexes with $C^{1}$-triangulations with $\dim N- \dim M$ odd.
\item A smooth map $f:M\rightarrow N$
\end{itemize}

{\bf Output}: `YES' if there is an immersion homotopic to $f$, `NO' otherwise.

{\bf Steps}:
\begin{enumerate}
  \item Using the algorithms in section \ref{effectiverep}, compute simplicial approximations of the classifying maps $\kappa_{M}:M\rightarrow BSO(M)$ and $\kappa_{N}:N\rightarrow BSO(N)$ for the corresponding tangent bundles.
  \item Construct the map $\phi:M\rightarrow BSO(M)\times BSO(N)$ where $\phi=(\kappa_{M}\times \kappa_{N})\circ (\id\times f)$.
  \item Construct $M^{+}$ and pick an extension $\phi^{+}$ of $\phi$ to $M^{+}$ (we note here that we have to include this step in general not only in the case that $M$ is simply connected because it is not in general decidable if $M$ is simply connected.)
  \item Using the relative minimal model for the appropriate codimension, use the algorithm from theorem \ref{liftingalg} to decide if the map $\phi^{+}$ lifts to $Mono(m-planes, n-planes)$.
  \item Output the result of the algorithm from the previous step.
\end{enumerate}

\printbibliography

\end{document}